\documentclass[12pt]{elsarticle}
\pdfoutput=1
\usepackage{booktabs,graphicx,multirow}
\usepackage{amsmath,amscd,amssymb,amsthm}

\makeatletter
\def\ps@pprintTitle{%
 \let\@oddhead\@empty
 \let\@evenhead\@empty
 \def\@oddfoot{}%
 \let\@evenfoot\@oddfoot}
\makeatother

\newtheoremstyle{standard}
  {3pt}           
  {3pt}           
  {\itshape}      
  {\parindent}    
  {}              
  {.}             
  {.5em}          
  {}              

\providecommand{\e}[1]{\ensuremath{\times 10^{#1}}}

\theoremstyle{standard}

\newtheorem{lemma}{LEMMA}[section]

\newtheorem{remark}{Remark}[section]

\newtheorem{observation}{Observation}[section]

\def\XXint#1#2#3{{\setbox0=\hbox{$#1{#2#3}{\int}$}
     \vcenter{\hbox{$#2#3$}}\kern-.5\wd0}}

\begin{document}


\begin{frontmatter}

\title
{
On the asymptotics of Bessel functions in the  Fresnel regime
}

\begin{abstract}
We introduce a version of the asymptotic expansions for Bessel functions
$J_\nu(z)$, $Y_\nu(z)$ that is valid whenever $|z| > \nu$ 
(which is deep in the Fresnel regime), 
as opposed to the standard expansions that 
are applicable only in the Fraunhofer regime (i.e. when $|z| > \nu^2$).
As expected, in the Fraunhofer regime our asymptotics reduce to the
classical ones. 
The approach is based on the observation that Bessel's equation admits a 
non-oscillatory phase function, and uses classical 
formulae to obtain an asymptotic expansion for this function;
this in turn leads to both an analytical tool and a numerical scheme for the 
efficient evaluation of $J_\nu(z)$, $Y_\nu(z)$, as well as 
various related quantities.
The effectiveness of the technique is demonstrated via several
numerical examples. We also observe that the procedure admits far-reaching
generalizations to wide classes of second order differential equations,
to be reported at a later date. 
\end{abstract}


\author[vr]{Zhu Heitman}
\author[jb]{James Bremer}
\author[vr]{Vladimir Rokhlin}
\author[bv]{Bogdan Vioreanu}


\address[vr]{Department of Computer Science, Yale University}
\address[jb]{Department of Mathematics, University of California, Davis}
\address[bv]{Department of Mathematics, University of Michigan}

\end{frontmatter}

Given a differential equation
\begin{equation}
y''(z) + q(z) y(z) = 0
\label{phase:eq1}
\end{equation}
on a (possible infinite) interval $(a,b)$,
a sufficiently smooth $\alpha$ is referred to as 
a phase function for~(\ref{phase:eq1}) if 
$\alpha'(z) \neq 0$ for all $z \in (a,b)$
and the pair of functions 
$u,v$ defined by the formulae
\begin{equation}
u(z) = \frac{\cos(\alpha(z)) }{\left|\alpha'(z)\right|^{1/2}},
\label{intro:u}
\end{equation}
\begin{equation}
v(z) = \frac{\sin(\alpha(z)) }{\left|\alpha'(z)\right|^{1/2}}
\label{intro:v}
\end{equation}
forms a basis in the  space of solutions  of~(\ref{phase:eq1}). %
Phase functions
arise from the theory of global transformations of ordinary
differential equations, which was initiated in \cite{Kummer};  more modern
discussions can be found, {\it inter alia}, 
in \cite{Boruvka}, \cite{Neuman}. Despite
their long history, phase functions possess a property that appears 
to have been overlooked: as long as the function $q$ 
is non-oscillatory, the equation~(\ref{phase:eq1}) possesses a
non-oscillatory phase function. This observation in its full generality
is somewhat technical, and will be reported at a later date; in this
short note, we apply it to the case of Bessel's equation 
\begin{equation}
z^2 \phi''(z) + z \phi'(z) + (z^2-\nu^2) \phi(z) = 0
\ \ \ \mbox{for all}\ \ 0 < z < \infty.
\label{bessel0}
\end{equation}

To the authors' knowledge,
the results in this paper are new despite the long-standing 
interest in the phase
functions associated with Bessel's equation. 
The most closely related antecedents of this work appear to 
be~\cite{Goldstein-Thaler}
and \cite{Spigler-Vianello2,Spigler-Vianello,Olver2}.
The former uses Taylor expansions of a non-oscillatory phase
function for Bessel's equation to evaluate the Bessel functions 
of large arguments; 
and the latter three works make use of phase functions 
to compute zeros of Bessel functions (as well as certain other 
special functions).

\begin{section}{Phase functions and the Kummer Equation}
If $u$ and $v$ are solutions of (\ref{phase:eq1}) related to 
a $C^3$ function $\alpha$ via the formulae
(\ref{intro:u}), (\ref{intro:v}),
then 
\begin{equation}
\tan(\alpha(z)) = \frac{v(z)}{u(z)}
\label{phase:1000}
\end{equation}
for all $z \in (a,b)$ such that $u(z) \neq 0$.   We differentiate (\ref{phase:1000}) 
in order to conclude that
\begin{equation}
\alpha'(z) \sec^2\left(\alpha(z)\right) = \frac{u(z)v'(z) - u'(z)v(z)}{u^2(z)}
\label{phase:1001}
\end{equation}
for all $z \in (a,b)$ such that $u(z) \neq 0$.
Note that the set of zeros of any nontrivial solution of (\ref{phase:eq1}) is discrete,
so that 
\begin{equation}
\{z \in (a,b) : u(z) \neq 0  \}
\end{equation}
is an open set.  From (\ref{intro:u}), (\ref{intro:v})
and a straightforward calculation we conclude that
\begin{equation}
u(z) v'(z)  - u'(z)v(z) = 1
\label{phase:1002}
\end{equation}
for all $z \in (a,b)$.   Moreover, by applying basic trigonometric identities
to (\ref{phase:1000}) we obtain the formula
\begin{equation}
\sec^2\left(\alpha(z)\right) =  \frac{u^2(z) + v^2(z)}{u^2(z)},
\label{phase:1003}
\end{equation}
which holds for all $z \in (a,b)$ such that $u(z) \neq 0$.
By combining (\ref{phase:1001}), (\ref{phase:1002}) and (\ref{phase:1003})
we conclude that
\begin{equation}
\alpha'(z) = \frac{1}{u^2(z) + v^2(z)}
\label{phase:alpha_prime}
\end{equation}
for all $z \in (a,b)$ for which $u(z)\neq 0$.  If $z \in (a,b)$ such
that $u(z) = 0$ but $v(z) \neq 0$, then 
\begin{equation}
\cot(\alpha(z)) = \frac{u(z)}{v(z)}.
\label{phase:cot}
\end{equation}
Using an argument analogous to that used to 
derive (\ref{phase:alpha_prime}) from formula (\ref{phase:1000}),
we conclude from (\ref{phase:cot}) that
(\ref{phase:alpha_prime}) also holds for all $z \in (a,b)$ such that $v(z) \neq 0$.
Since (\ref{phase:1002}) implies that $u$ and $v$ are never simultaneously
zero on the interval $(a,b)$, formula (\ref{phase:alpha_prime}) in fact holds
for all $z \in (a,b)$.

By differentiating (\ref{intro:u}) twice we obtain
\begin{equation}
u''(z) = 
\frac{\cos(\alpha(z))}{\sqrt{\left|\alpha'(z)\right|}}
\left(
- (\alpha'(z))^2
- \frac{1}{2}\left(\frac{\alpha'''(z)}{\alpha'(z)}\right)
+ \frac{3}{4}\left(\frac{\alpha''(z)}{\alpha'(z)}\right)^2
\right).
\label{phase:3000}
\end{equation}
By adding $q(z) u(z)$ to both sides of (\ref{phase:3000})
and making use of (\ref{intro:u}) we conclude that
%
\begin{equation}
u''(z) +q(z) u(z)= 
\frac{\cos(\alpha(z))}{\sqrt{\left|\alpha'(z)\right|}}
\left(
q(z) 
- (\alpha'(z))^2
- \frac{1}{2}\left(\frac{\alpha'''(z)}{\alpha'(z)}\right)
+ \frac{3}{4}
\left(\frac{\alpha''(z)}{\alpha'(z)}\right)^2
\right).
\label{phase:3001}
\end{equation}
Using an analogous sequence of steps, we obtain the formula
\begin{equation}
v''(z) +q(z) v(z)= 
\frac{\sin(\alpha(z))}{\sqrt{\left|\alpha'(z)\right|}}
\left(
q(z) 
- (\alpha'(z))^2
- \frac{1}{2}\left(\frac{\alpha'''(z)}{\alpha'(z)}\right)
+ \frac{3}{4}
\left(\frac{\alpha''(z)}{\alpha'(z)}\right)^2
\right).
\label{phase:3002}
\end{equation}
from  (\ref{intro:v}).
Suppose that  $\alpha$ is a phase function for (\ref{phase:eq1}).
Then
\begin{equation}
u''(z) + q(z) u(z) = 0 = v''(z) + q(z) v(z) 
\label{phase:3003}
\end{equation}
for all $z \in (a,b)$.  Since 
\begin{equation}
\cos^2(\alpha(z)) + \sin^2(\alpha(z)) = 1,
\end{equation}
the functions
\begin{equation}
\cos(\alpha(z))
\end{equation}
and
\begin{equation}
\sin(\alpha(z))
\end{equation}
cannot simultaneously be $0$.  By combining this observation with 
(\ref{phase:3001}), (\ref{phase:3002}) and (\ref{phase:3003}),
we conclude that  $\alpha$ satisfies  the third order nonlinear differential equation

\begin{equation}
\left(\alpha'(z)\right)^2 = q(z) - \frac{1}{2} \frac{\alpha'''(z)}{\alpha'(z)} + 
\frac{3}{4} \left(\frac{\alpha''(z)}{\alpha'(z)}\right)^2
\label{phase:kummer}
\end{equation}
on the interval $(a,b)$.  Suppose, on the other hand, that 
 $\alpha$ is an element
of $C^3\left(\left(a,b\right)\right)$ such that $\alpha'(z) \neq 0$
for all $z \in (a,b)$, and that
the function $q$ is defined by equation (\ref{phase:kummer}).    Then 
 (\ref{phase:3001})
and (\ref{phase:3002}) imply that the functions $u$ and $v$ 
defined via  (\ref{intro:u}),
(\ref{intro:v})
are solutions of (\ref{phase:eq1}).  A straightforward calculation
shows that the Wronskian 
\begin{equation}
u(z)v'(z)-u(z)v(z) 
\end{equation}
of $\{u,v\}$   is $1$, so that the pair of functions $u$, $v$ in fact
forms a basis in the 
space of solutions
of (\ref{phase:eq1}).
We will refer to 
(\ref{phase:kummer})
as Kummer's equation, after E.~E. Kummer who studied it in~\cite{Kummer}.

Our principal interest is in the highly oscillatory case,
where $q(z)$ is of the form $\gamma^2 \widetilde{q}(z)$ with $\gamma$ 
a large real-valued constant
and $\widetilde{q}$ a complex-valued function, so that 
$q(z)$ is asymptotically of the order $\gamma^2$.
As a consequence of the Sturm 
comparison theorem (see, for example,~\cite{Coddington-Levinson})
 solutions of (\ref{phase:eq1}) are necessarily highly oscillatory
when $\gamma$ is large.   
Moreover, the form of (\ref{phase:kummer}) and the appearance
of $\lambda$ in it suggests 
that  phase functions for (\ref{phase:eq1}) 
will also be highly oscillatory in this regime  --- and most of them are.
However, non-oscillatory $\alpha$ are
available whenever $q$ is non-oscillatory (in an appropriate sense). 
In this note, we construct such an $\alpha$ in the case of Bessel's
equation.  

%


\label{section:phase}
\end{section}


\begin{section}{A Phase function for Bessel's equation}\label{phase_section}
The Bessel function of the first kind of order $\nu$
\begin{equation*}
J_\nu(z) = 
\left(\frac{z}{2}\right)^\nu
\sum_{n=0}^\infty \frac{(-\frac{1}{4}z^2)^n}{n!\ \Gamma\left(\nu+n+1\right)} 
\end{equation*}
and the Bessel function of the second kind of order $\nu$
\begin{equation*}
Y_\nu(z) = \frac{J_\nu(z) \cos(\nu \pi) - J_{-\nu}(z)}{\sin(\nu \pi)}
\end{equation*}
form a basis in the space of solutions of Bessel's equation~(\ref{bessel0}).

Bessel functions are among the most studied and well-understood of 
special functions. Even the standard reference books such 
as~\cite{Abramowitz-Stegun},~\cite{Gradstein} contain a wealth of 
information. Here, we only list facts necessary for the specific 
purposes of the paper; the reader is referred 
to~\cite{Watson},~\cite{Abramowitz-Stegun},~\cite{Gradstein} (and 
many other excellent sources) for more detailed information.

For arbitrary positive real $\nu$ and large values of 
$|z|$, 
the functions $J_{\nu}$, $Y_{\nu}$ possess asymptotic expansions 
of the form
\begin{equation}
J_{\nu} (z) = \sqrt{{2 \over \pi \, z}} \, \cos(z-{\nu \pi \over 2}
            -{\pi \over 4}) \, \left(1+\sum_{k=1}^{\infty} 
              {P_{2 \, k} \, (\nu) \over z^k} \right),
\label{asympj}
\end{equation}
\begin{equation}
Y_{\nu} (z) = \sqrt{{2 \over \pi \, z}} \, \sin(z-{\nu \pi \over 2}
            -{\pi \over 4}) \, \left(1+\sum_{k=1}^{\infty} 
              {Q_{2 \, k} \, (\nu) \over z^k} \right),
\label{asympy}
\end{equation}
where for each $k=1,2,\cdots$, $P_{2k},Q_{2k}$ are polynomials of 
order $2k$ (see,
for example,~\cite{Abramowitz-Stegun} for exact definitions
of $\{P_{2k} \}$, $\{Q_{2k} \}$).

\begin{remark}

As one would expect from the form of
expansions~(\ref{asympj}),~(\ref{asympy}), they provide no useful
approximations to $J_{\nu}(z)$, $Y_{\nu}(z)$ 
until $|z|$ is greater than roughly $\nu^2$.
The regime where $\nu < |z| < \nu^2$ is known as {\it Fresnel regime},
while the regime $|z| > \nu^2$ is referred to as 
{\it Fraunhofer regime}. In other words, in the Fraunhofer regime the 
classical asymptotics~(\ref{asympj}),~(\ref{asympy}) become
useful.

\end{remark}

The transformation $\psi(z) = z^{1/2}\varphi(z)$  brings~(\ref{bessel0})
into the standard form
\begin{equation}
\psi''(z) + \left(1-\frac{\nu^2-1/4}{z^2} \right) \psi(z) = 0.
\label{bessel:standard}
\end{equation}
Two standard solutions of~(\ref{bessel:standard}) are 
$\{\sqrt{z}J_\nu(z),\sqrt{z}Y_\nu(z)\}$; in this 
case,~(\ref{phase:alpha_prime}) becomes
\begin{equation}
\alpha_\nu'(z) =  \frac { 2 }{\pi z}\frac{1}{J_\nu^2(z) + Y_\nu^2(z)}
\label{bessel:bessel_prime}
\end{equation}
Clearly,~(\ref{bessel:bessel_prime}) determines the phase
function $\alpha$ up to a constant; choosing 
\begin{equation}
\alpha_\nu(z) = -\frac{\pi}{2} + \int_0^z \alpha_\nu'(u) du,
\label{bessel:particular_alpha}
\end{equation}
we obtain expressions 
\begin{equation}
J_\nu(z)= M_\nu(z) \cos(\alpha_\nu(z)) 
\label{forJ}
\end{equation}
\begin{equation}
Y_\nu(z)=M_\nu(z) \sin(\alpha_\nu(z)), 
\label{forY}
\end{equation}
with $M_\nu(z)$ defined by the formula
\begin{equation}
M_\nu(z) = \sqrt{ J_\nu^2(z) + Y_\nu^2(z) }.
\label{M_first}
\end{equation}

In a remarkable coincidence, there exists a simple integral
expression for $M_{\nu}$, valid for all z such that $\arg(z) < \pi$.
Specifically, 
\begin{equation}
(M_{\nu} (z))^2 = J_\nu^2(z) + Y_\nu^2(z) = \frac{8}{\pi^2}
\int_0^\infty K_0(2z\sinh(t)) \cosh(2\nu t)\ dt
\label{nickol}
\end{equation}
(see, for example,~\cite{Gradstein}, Section 6.664); even a cursory examination 
of~(\ref{nickol}) shows that for $z$ on the real axis, 
$M_{\nu}$ is a non-oscillatory function of $z$.

The approximation 
\begin{equation}
M_\nu^2(z) \sim \frac{2}{\pi z} 
\sum_{n=0}^\infty 
\frac{\Gamma\left(n+\frac{1}{2}\right)}{n!\sqrt{\pi}} 
\frac{\Gamma\left(\nu +\frac{1}{2}+n\right)}
{\Gamma\left(\nu+\frac{1}{2}-n\right)}
\frac{1}{z^{2n}},
\label{bessel:amp}
\end{equation}
can be found (for example) in Section~13.75 of~\cite{Watson}; 
it is asymptotic in ${1 \over z}$, and its first several terms are
\begin{equation}
 M_\nu^2(z) \sim 
\frac{2}{\pi z}
\left( 1 + 
\frac{1}{2} \frac{\mu - 1}{(2z)^2} +
\frac{1}{2}\cdot \frac{3}{4} \frac{(\mu - 1)(\mu-9)}{(2z)^4} +
\frac{1}{2}\cdot \frac{3}{4}\cdot \frac{5}{6} 
\frac{(\mu-1)(\mu-9)(\mu-25)}{(2z)^6} +\cdots\right),
\label{bessel:amplitude}
\end{equation}
with $\mu = 4\nu^2$. When the expansion~(\ref{bessel:amp})
is truncated after $k$ terms, the error of the resulting approximation
is bounded by the absolute value of the next term, as long as 
$z$ is real and $k > \nu$ (see~\cite{Watson}). Computationally, it
is often convenient to rewrite~(\ref{bessel:amp}) in the form 
\begin{equation}
M_\nu^2(z) \sim \frac{2}{\pi z}\left( 1+\sum_{n=1}^\infty 
{t_n \over z^{2n}} \right),
\label{bessel:amp2}
\end{equation}
with $t_0=1$, and 
\begin{equation}
t_{n} = t_{n-1} \left(\frac{\mu - (2n-1)^2}{4}\right) \frac{2n-1}{2n}.
\label{t_n0}
\end{equation}

\begin{observation}

While the expansion~(\ref{bessel:amplitude}) has been known for
almost a century, one of its implications does not appear to be
widely understood. Specifically,~(\ref{bessel:amplitude}) provides
an effective approximation to  $M_\nu$ whenever $|z| > \nu$;
combined with~(\ref{forJ}),~(\ref{forY}), it results in asymptotic
approximation to $J_{\nu}$ ,$Y_{\nu}$  in 
the Fresnel regime, i.e. when $|z| > \nu$
- as opposed to the expansions~(\ref{asympj}),~(\ref{asympy}), 
that are only valid in the Fraunhofer regime
- i.e. when $|z| > \nu^2$. This (rather elementary) observation
has obvious implications, inter alia,  for the numerical evaluation of Bessel
functions of high order.

\end{observation}

Applying Lemma~\ref{lemma1} to the expansion~(\ref{bessel:amp2})
and using~(\ref{bessel:bessel_prime}), (\ref{M_first}),
we obtain 
\begin{equation}
\alpha'_\nu(z) =  \frac{2}{\pi z} {1 \over M_\nu^2(z)} 
\sim  \left( 1+\sum_{n=1}^\infty 
{s_n \over z^{2n}} \right),
\label{bessel:amp3}
\end{equation}
with $s_0=1$, and $s_1,s_2,s_3, \ldots$ defined by the formula
\begin{equation}
s_n = - \left( s_n + \sum_{j=1}^{n-1} t_j s_{n-j} \right);
\label{s_n0}
\end{equation}
the first several terms in~(\ref{s_n0}) are 
\begin{equation}
\alpha'_\nu(z) \sim
1 -  \frac{\mu-1}{8z^2} - \frac{\mu^2-26\mu+25}{128 z^4} 
-\frac{\mu^3-115 \mu^2 + 1187 \mu - 1073}{1024 z^6}
+\cdots.
\label{bessel:alpha_asymptotic_der}
\end{equation}
Obviously, the indefinite integral of~(\ref{bessel:amp3}) is
\begin{equation}
\alpha_\nu(z) 
\sim C+z-\sum_{n=1}^\infty 
{s_n \over (2 n -1) \, z^{2n-1}},
\label{bessel:integr1}
\end{equation}
with $C$ to be determined; the first several terms in~(\ref{bessel:integr1}) 
are 
\begin{equation}
\alpha_\nu(z) \sim
C+z +  \frac{\mu-1}{8z} + \frac{\mu^2-26\mu+25}{384 z^3} 
+\frac{\mu^3-115 \mu^2 + 1187 \mu - 1073}{5120 z^5}
+\cdots.
\label{bessel:alpha_asymptotic}
\end{equation}
In order to find the value of $C$ 
in~(\ref{bessel:integr1}),~(\ref{bessel:alpha_asymptotic}),
we observe that for sufficiently
large $|z|$, the expansion~(\ref{bessel:alpha_asymptotic}) becomes
\begin{equation}
\alpha_\nu(z) \sim C+z.
\label{simpler}
\end{equation}
Substituting~(\ref{simpler}) into~(\ref{forJ}),~(\ref{forY}), we 
have 
\begin{equation}
J_{\nu} (z) \sim M_\nu(z) \cos(C+z),
\label{simplerj}
\end{equation}
\begin{equation}
Y_{\nu} (z) \sim M_\nu(z) \sin(C+z).
\label{simplery}
\end{equation}
Clearly, the approximations~(\ref{simplerj}),~(\ref{simplery}) 
must be compatible with~(\ref{asympj}),~(\ref{asympy}), which
yields 
\begin{equation}
C= -{\nu \pi \over 2}-{\pi \over 4}. 
\label{cvalue}
\end{equation}
Finally, substituting~(\ref{cvalue}) into~(\ref{bessel:alpha_asymptotic}),
we end up with 
\begin{equation}
\alpha_\nu(z) 
\sim -{\nu \pi \over 2}-{\pi \over 4} +z-\sum_{n=1}^\infty 
{s_n \over (2 n -1) \, z^{2n-1}},
\label{bessel:integr2}
\end{equation}
or
\begin{equation}
\alpha_\nu(z) \sim
-{\nu \pi \over 2}-{\pi \over 4}
+z +  \frac{\mu-1}{8z} + \frac{\mu^2-26\mu+25}{384 z^3} 
+\frac{\mu^3-115 \mu^2 + 1187 \mu - 1073}{5120 z^5}
+\cdots,
\label{alpha_asymp}
\end{equation}
with $\mu=4 \, \nu^2$. We note that the asymptotic expansion
(\ref{bessel:integr2}) is only useful because
the particular phase function $\alpha_\nu(z)$ it represents is non-oscillatory:
Figure~\ref{figure:phase} shows plots of the derivative
of the phase function associated with $\{\sqrt{z} J_\nu(z),\sqrt{z} Y_\nu(z)\}$
and also the derivative of the phase function 
associated with another choice of basis, 
$\{2 \sqrt{z} J_\nu(z),\sqrt{z} Y_\nu(z)\}$.


\label{section:bessel}
\end{section}

\begin{figure}[b!!]

\begin{center}

\includegraphics[width=.70\textwidth]{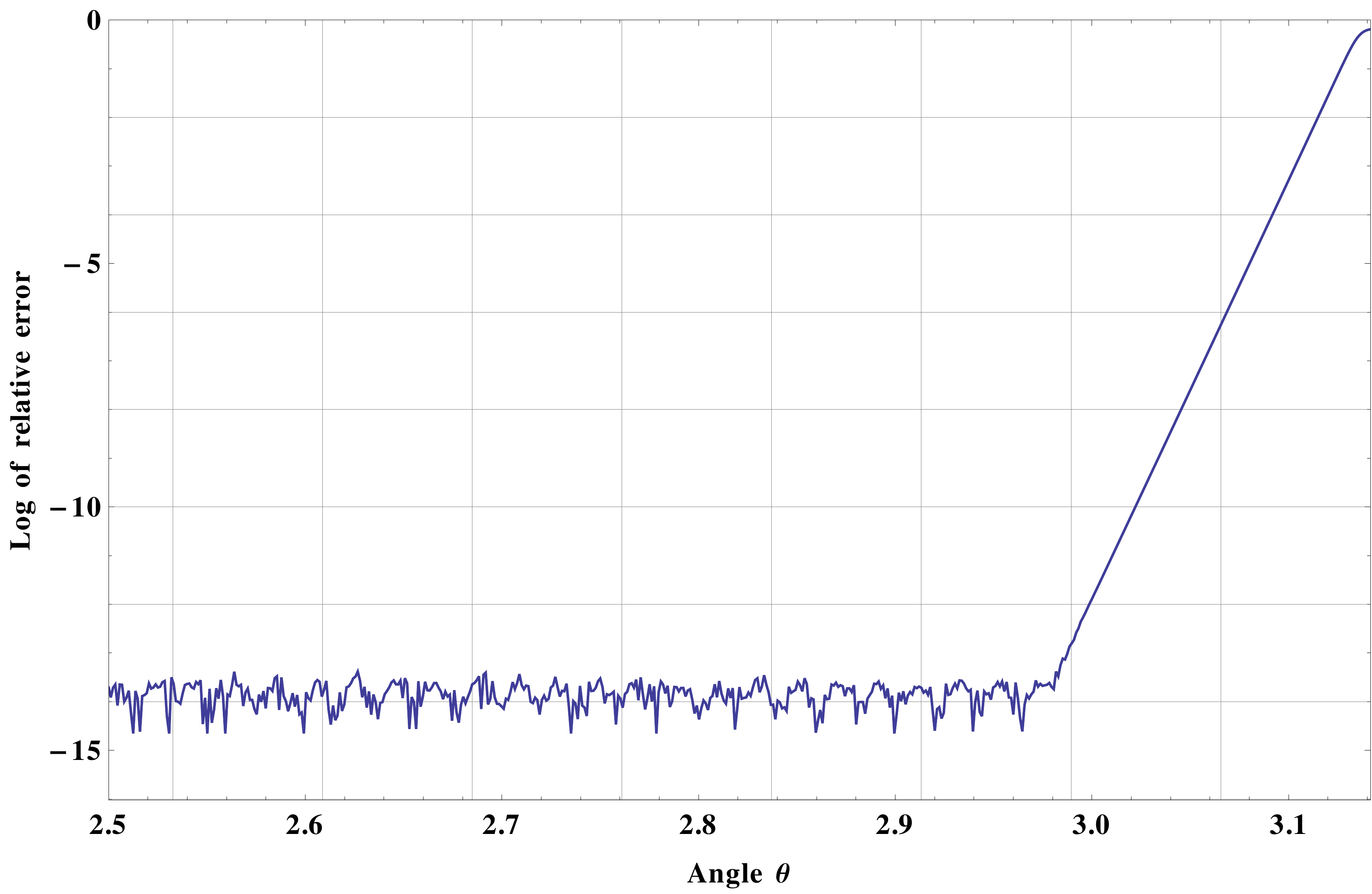}
\caption{{\bf Failure as $\mathbf{\mbox{\bf arg}(z)}$ approaches $\mathbf{\pi}$}.  
The base-10 logarithm of the 
 relative error in the approximation of $Y_{10}(100\exp(i\theta))$
as $\theta$ approaches $\pi$.}
\label{figure:logplot2}

\end{center}

\end{figure}

\begin{figure}[h!!]

\begin{center}

\includegraphics[width=.7\textwidth]{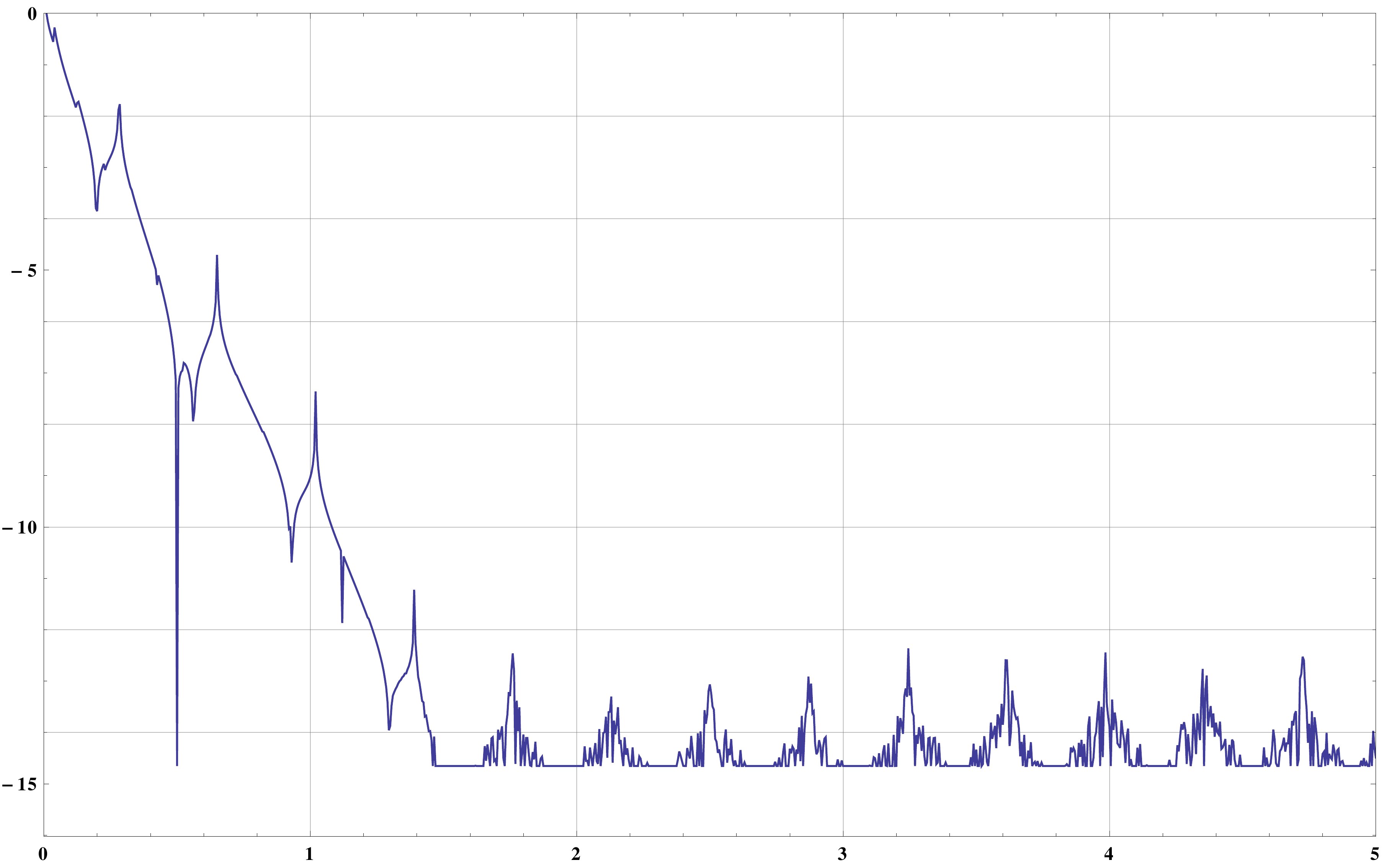}
\caption{{\bf Failure for small orders}.  The base-10 logarithm of the 
 relative error in the approximation of $J_\nu(10\nu)$ 
plotted as a function of $\nu$.}
\label{figure:logplot}

\end{center}

\end{figure}

\begin{section}{Numerical Experiments}

In this section we present the result of several numerical experiments
conducted to verify the scheme of this paper.  The code for
these experiments was written in Fortran~77 and compiled using
the Intel Fortran compiler version 12.0.  Experiments were
conducted on a laptop equipped with an Intel Core i7-2620M processor
running at 2.70 GHz and 8 GB of RAM.  Machine zero was
$\varepsilon_0 = 2.22044604925031\e{-16}$.

We constructed asymptotic approximations to the functions $J_{\nu},
Y_{\nu}$ via the formulae~(\ref{t_n0}),~(\ref{s_n0}),~(\ref{bessel:integr2}).
In order to avoid exceeding the machine exponent, we altered the procedure 
slightly, so that the coefficients $t_k, s_k$ are never computed by themselves:
only the ratios 
\begin{equation}
{t_k \over z^{2k}}, \ {s_k \over z^{2k}}
\label{rats}
\end{equation}
are calculated via obvious modifications of~(\ref{t_n0}),~(\ref{s_n0}).

As with all procedures relying on asymptotic expansions,
it is not always possible to achieve a desired accuracy.
Indeed, the magnitudes of the terms of the expansion 
reach a certain minimum 
and then proceed to increase.  And, of course, truncating the expansions when
the terms become small does not necessarily
ensure the accuracy of the approximation (see~\cite{Olver} for 
numerous examples of possible pathologies).
It is shown in Section~13.75 of \cite{Watson} that if $\nu$ is real, 
$z$ is positive
and $n > \nu-1/2$, then the remainder resulting
from the first $n$ terms of  expansion
(\ref{bessel:amplitude})
 is smaller in magnitude than the $(n+1)$st term.  However, the authors
are not aware of any error bounds for (\ref{bessel:integr2})
(or for~(\ref{bessel:amplitude}) ) in the general case. 


Moreover, the value of the phase function $\alpha_\nu(z)$ is proportional
to the argument $z$ and values of $J_\nu(z)$ and $Y_\nu(z)$ are obtained
in part by evaluating the sine and cosine of $\alpha_\nu(z)$.  This imposes
limitations on the accuracy of the obtained approximations
when $z$ is large due to the well-known
difficulties in evaluating periodic functions of large arguments.

\begin{subsection}{Comparison with Mathematica}

In these first experiments, we applied the procedure of this paper 
to the evaluation of 
the Bessel functions $J_\nu(z)$ and $Y_\nu(z)$ at various
orders and arguments.  The resulting values were compared with those produced
by version 9.0.0 of Wolfram's Mathematica package; 30 digit precision 
was requested from Mathematica. 
Table~\ref{table:experiments_one} reports the results.  There,
the number of terms used in the expansions of the modulus and phase
functions and the relative errors in the obtained
values of $J_\nu(z)$ and $Y_\nu(z)$ are reported. 

\label{section:experiments_mathematica}

\end{subsection}

\begin{subsection}{Bessel functions of large order.}
In this experiment, we 
approximate the values of Bessel functions of very large orders
and arguments.  Comparison with other approaches is difficult for such
large orders; for instance, Mathematica's Bessel function routines
are prohibitively slow in this regime.
We settled for running our procedure twice, once using double precision
arithmetic and once using extended precision (Fortran REAL*16) arithmetic
in order to produce reference values for comparison.  

The first row of each entry in Table~\ref{table:large_orders1} reports the relative error in 
the approximations
of $J_\nu(z)$ and the second row gives the relative error in the approximation of $Y_\nu(z)$.
 Table~\ref{table:large_orders2} 
gives the number of terms in the expansion
of the modulus function and the number of terms in the expansions of the modulus and
phase function used to evaluate $J_\nu(z)$ and $Y_\nu(z)$.
These values depended only on the ratio of $ z$ to $\nu$ and not on the value of $\nu$. 

\end{subsection}

\begin{subsection}{Failure for small orders.}

In this experiment, we considered the performance of the procedure of
this paper for relatively small values of  $|\nu|$.
We evaluated $J_\nu(10\nu)$ at a series of values 
of $\nu$ between $0$ and $5$.  A plot
of the base-$10$ logarithm of the 
relative error in $J_\nu(10\nu)$ is shown
in Figure~\ref{figure:logplot}. 
Errors were estimated via comparison with Wolfram's Mathematica
package; 30 digit precision was requested from Mathematica.

\end{subsection}

\begin{subsection}{Failure as $\mathbf{\mbox{\bf arg}(z)}$ approaches $\mathbf{\pi}$.}

In this experiment, we computed the values of $Y_{10}(\exp(i \theta))$ 
as $\theta $ approaches $\pi$.
The obtained values were compared to those
reported by Wolfram's Mathematica package; 30 digit precision
was once again requested from Mathematica.  A plot of the base-$10$ logarithm of the 
 relative 
error in $Y_{10}(\exp(i\theta))$ as a function of 
$\theta$ is shown in Figure~\ref{figure:logplot2}. 

\end{subsection}
\label{section:experiments}
\end{section}

\begin{section}{Conclusions}
We have shown that the Bessel functions $J_\nu(z)$ and $Y_\nu(z)$ 
can be efficiently evaluated when $\nu$ is large and $|z| > |\nu|$.
This was achieved by representing
Bessel functions in terms of a non-oscillatory phase function for which
(conveniently enough) a well-known asymptotic expression is available.

The observation underlying the scheme of this paper --- namely, 
the existence of a
non-oscillatory phase function --- is not a peculiarity of Bessel's equation.
The solutions of a large class of second order linear differential equations
can be approximated to high accuracy via non-oscillatory phase functions,
a development the authors will report at a later date.
\label{section:conclusion}
\end{section}

\begin{section}{Appendix}
Here, we formulate a lemma used in Section~\ref{phase_section}; its proof is 
an exercise in elementary calculus, and can be found, for example, 
in~\cite{Olver}.

\begin{lemma}
Suppose that
\begin{equation}
f(z) \sim 1 + \frac{a_1}{z^2} + \frac{a_2}{z^4} + \frac{a_3}{z^6} + \cdots
\end{equation}
is an asymptotic expansion for $f: C \to C$, with $a_1,a_2,a_3,\ldots$ 
a sequence of complex numbers. Then the asymptotic expansion of \, $1/f$
is 
\begin{equation}
\frac{1}{f(z)} \sim 1 + \frac{b_1}{z^2} + \frac{b_2}{z^4} + 
\frac{b_3}{z^6} + \cdots,
\end{equation}
where $b_1=-a_1$, and the rest of the coefficients $b_2,b_3,b_4,\ldots$ 
are given by the formula 
\begin{equation}
b_n = -a_n -\sum_{j=1}^{n-1} a_{j} b_{n-j}.
\end{equation}
\label{lemma1}
\end{lemma}

\end{section}
   
\begin{section}{Acknowledgments}
We would like to thank the reviewer for a careful reading of the manuscript
and for several useful suggestions.
Zhu Heitman was supported in part by the Office of Naval Research
under contracts ONR N00014-10-1-0570 and ONR N00014-11-1-0718.
James Bremer was supported in part by a fellowship from the Alfred P. Sloan
Foundation and by National Science Foundation
 grant DMS-1418723.
  Vladimir Rokhlin was supported in part by 
Office of Naval Research contracts
ONR N00014-10-1-0570 and ONR N00014-11-1-0718, and by the Air Force Office of Scientific
Research under contract AFOSR FA9550-09-1-0241.
\end{section}

\bibliographystyle{acm}
\bibliography{bessel}

\vfill 
\eject

\begin{figure}[t!!]
\includegraphics[width=.49\textwidth]{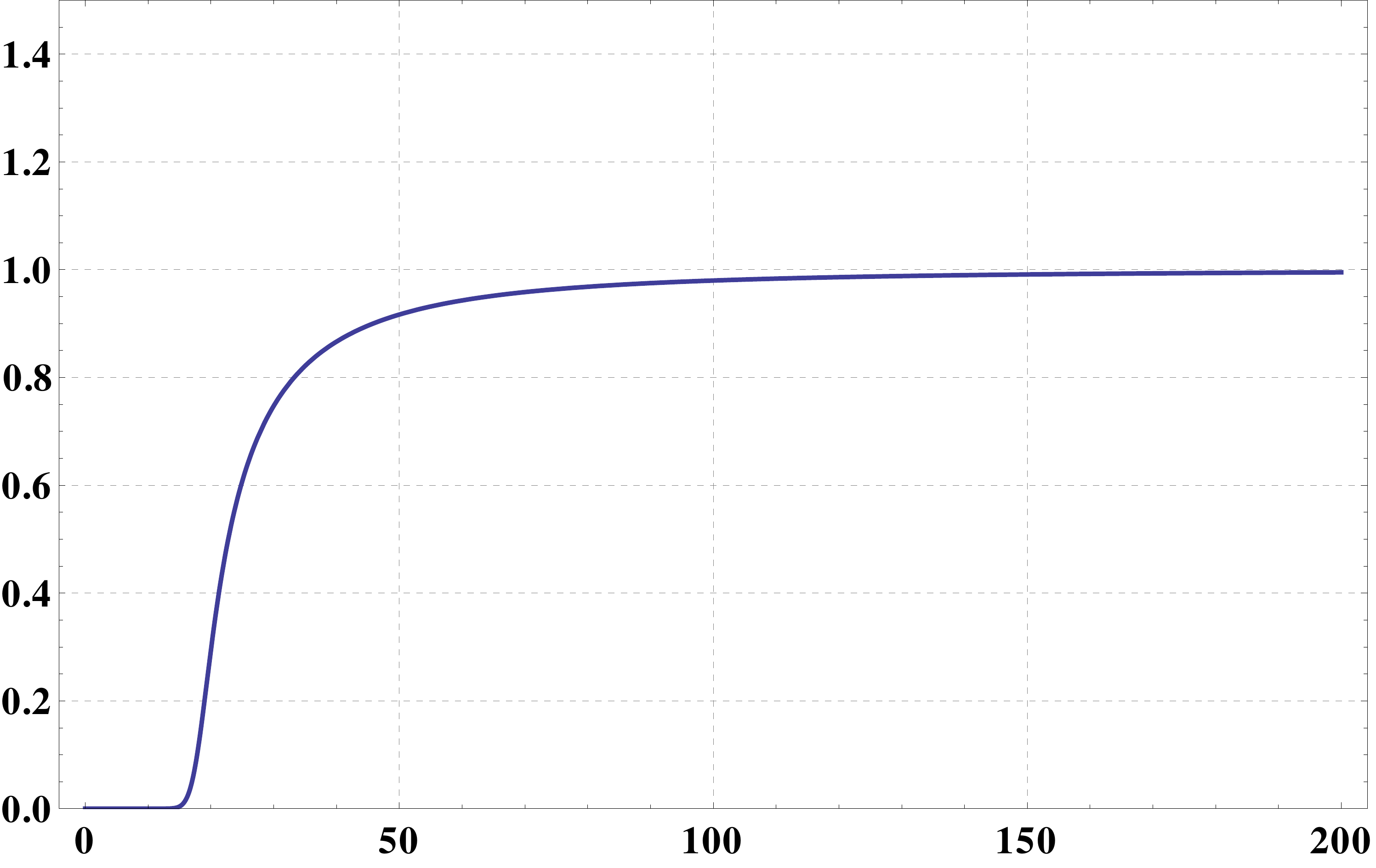}
\hfill
\includegraphics[width=.49\textwidth]{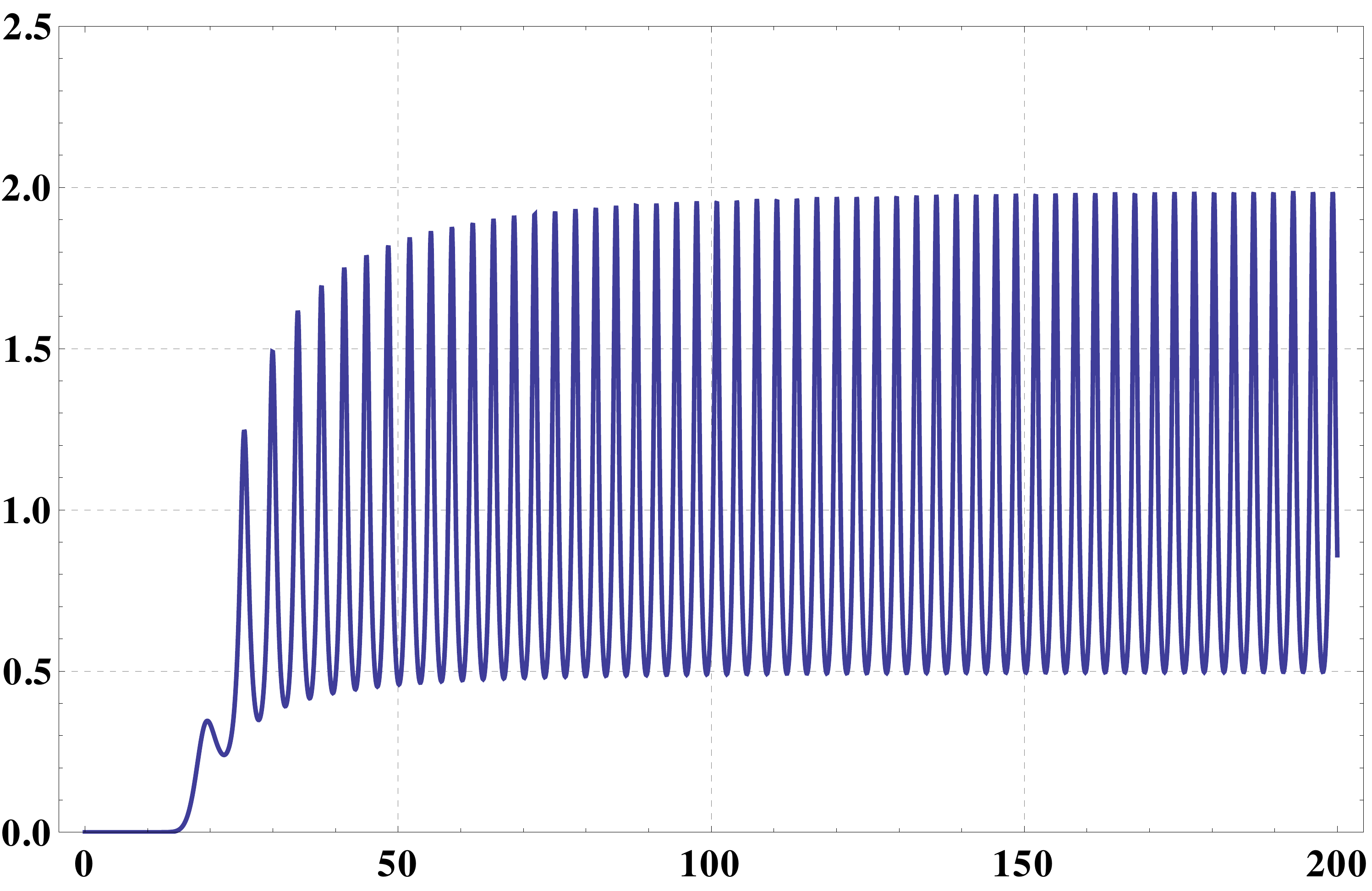}
\caption{On the left, a plot of the derivative 
of a phase function
associated with the basis of solutions 
$\{\sqrt{z} J_{20}(z), \sqrt{z} Y_{20}(z)\}$
for (\ref{bessel:standard}); 
 on the right, a plot of the derivative of a phase function associated 
with the basis  $\{2 \sqrt{z} J_{20}(z), \sqrt{z} Y_{20}(z)\}$.
Both are over the interval $[0,200]$.}
\label{figure:phase}
\end{figure}

\begin{table}[h!!!]

\vspace{1in}
\begin{center}

\begin{tabular}{ccc}
$z/\nu$  & Modulus terms & Phase terms \\
\midrule 
\addlinespace[.5em]
$1.1$    & $162$         & $133$ \\
$2$      & $024$         & $022$ \\
$10$     & $008$         & $008$ \\
$10\pi$  & $006$         & $006$ \\
\addlinespace[.5em]
\end{tabular}
\caption{{\bf Bessel functions of large order}.  The number of terms in 
the expansions of the modulus
and phase functions used in the computation of Bessel functions of large orders.}
\label{table:large_orders2}
\vskip -2em

\end{center}

\end{table}

\begin{table}[h!!!]

\begin{center}
\small
\begin{tabular}{rrcccc}
      &     & Modulus      & Phase &  Relative error & Relative error\\
$\nu$ & $z$ & terms        & terms &  in $J_\nu(z)$  & in $Y_\nu(z)$ \\
\midrule
\addlinespace[.5em]
$50$          & $1.1|\nu|$ & $049$  & $091$ &  $1.45\e{-14}$ & $3.25\e{-14}$ \\
              & $2|\nu|$   & $023$  & $022$ &  $3.73\e{-14}$ & $9.80\e{-15}$ \\
              & $10|\nu|$  & $008$  & $008$ &  $1.77\e{-14}$ & $9.38\e{-15}$ \\
              & $100|\nu|$ & $005$  & $005$ &  $4.53\e{-13}$ & $7.23\e{-14}$ \\
\addlinespace[.5em]
$50-10i$      & $1.1|\nu|$ & $075$  & $098$ &  $4.85\e{-15}$ & $4.81\e{-15}$ \\
              & $2|\nu|$   & $023$  & $022$ &  $6.04\e{-15}$ & $6.12\e{-15}$ \\
              & $10|\nu|$  & $008$  & $008$ &  $9.79\e{-14}$ & $9.77\e{-14}$ \\
              & $100|\nu|$ & $005$  & $005$ &  $6.31\e{-13}$ & $6.31\e{-13}$ \\
\addlinespace[.5em]
$100+20i$     & $1.1|\nu|$ & $080$  & $117$ &  $1.79\e{-14}$ & $1.79\e{-14}$ \\
              & $2|\nu|$   & $024$  & $021$ &  $3.19\e{-14}$ & $3.19\e{-14}$ \\
              & $10|\nu|$  & $008$  & $008$ &  $1.46\e{-13}$ & $1.46\e{-13}$ \\
              & $100|\nu|$ & $005$  & $005$ &  $5.43\e{-13}$ & $5.43\e{-13}$ \\
\addlinespace[.5em]
$10000$       & $1.1|\nu|$ & $162$  & $133$ &  $1.72\e{-11}$ & $4.63\e{-12}$ \\
              & $2|\nu|$   & $024$  & $022$ &  $1.18\e{-12}$ & $6.69\e{-13}$ \\
              & $10|\nu|$  & $008$  & $008$ &  $8.96\e{-13}$ & $1.28\e{-10}$ \\
              & $100|\nu|$ & $005$  & $005$ &  $1.10\e{-10}$ & $2.96\e{-12}$ \\
\addlinespace[.5em]
$100000$      & $1.1|\nu|$ & $162$  & $133$ &  $5.18\e{-13}$ & $7.16\e{-13}$ \\
              & $2|\nu|$   & $024$  & $022$ &  $4.82\e{-11}$ & $2.34\e{-11}$ \\
              & $10|\nu|$  & $008$  & $008$ &  $1.57\e{-09}$ & $3.93\e{-12}$ \\
              & $100|\nu|$ & $005$  & $005$ &  $1.92\e{-10}$ & $2.41\e{-10}$ \\

\end{tabular}
\vskip 1em
\caption{{\bf Comparison with Mathematica.}  The results of the experiments of Section~\ref{section:experiments_mathematica}.}
\label{table:experiments_one}

\end{center}
\end{table}

\begin{table}[h!!!]

\begin{center}

  \small
\begin{tabular}{cccccc}
$z$                        & $\nu = 10^6$  & $\nu = 10^9$   & $\nu = 10^{12}$ & $\nu = 10^{15}$ & $\nu = 10^{18}$ \\
\midrule 
\addlinespace[.5em]
\multirow{2}{*}{$1.1\nu$}  & $1.96\e{-11}$ & $4.48\e{-07}$  & $3.74\e{-07}$  & $3.26\e{-08}$   & $5.78\e{-06}$ \\
                           & $1.18\e{-10}$ & $2.81\e{-07}$  & $1.70\e{-08}$  & $1.57\e{-06}$   & $5.87\e{-08}$ \\
\addlinespace[.5em]
\multirow{2}{*}{$2\nu$}    & $1.04\e{-10}$ & $2.15\e{-07}$  & $1.01\e{-07}$  & $1.25\e{-07}$   & $4.08\e{-07}$ \\
                           & $4.66\e{-11}$ & $6.94\e{-07}$  & $1.37\e{-07}$  & $2.96\e{-07}$   & $7.73\e{-07}$ \\
\addlinespace[.5em]
\multirow{2}{*}{$10\nu$}   & $5.25\e{-09}$ & $3.80\e{-07}$  & $3.20\e{-05}$  & $1.09\e{-05}$   & $1.04\e{-06}$ \\
                           & $4.53\e{-10}$ & $5.30\e{-07}$  & $4.33\e{-09}$  & $6.83\e{-08}$   & $4.00\e{-06}$ \\
\addlinespace[.5em]
\multirow{2}{*}{$10\pi\nu$}& $3.23\e{-10}$ & $4.58\e{-06}$  & $7.33\e{-06}$  & $3.39\e{-06}$   & $4.12\e{-06}$ \\
                           & $3.80\e{-10}$ & $1.01\e{-05}$  & $2.10\e{-07}$  & $3.68\e{-06}$   & $3.05\e{-06}$ \\

\addlinespace[.5em]
\end{tabular}
\caption{{\bf Bessel functions of large order}.  The relative errors in the 
obtained approximations of Bessel functions of large orders: the round-off 
due to large values of $|z|$ sharply limits the accuracy of results 
obtainable in double precision}
\label{table:large_orders1}
\vskip -2em

\end{center}

\end{table}

\end{document}